\documentclass[3p,12pt]{elsarticle}

\usepackage{times,amsmath,amssymb,multirow,bm,graphicx,hyperref}
\usepackage[latin1]{inputenc}
\usepackage{colortbl,rotating,lscape}
\usepackage[english]{babel}

\newcommand{\deljump}{[\![\delta]\!]\>}

\begin{document}

\title{The effect of the elliptic polarization on the quasi-particle dynamics of
linearly coupled systems of Nonlinear Schr\"odinger Equations}

\author[mdt]{M. D. Todorov\corref{cor1}}
\ead{mtod@tu-sofia.bg}
\address[mdt]{Dept of Diff. Equations, Faculty of Applied
Math. and Informatics,  Technical University of Sofia, Sofia, 1000, Bulgaria}


\cortext[cor1]{Corresponding author. Tel.: +359-899-187089.}



 \begin{abstract}
We investigate numerically by a conservative difference scheme in complex arithmetic the head-on and takeover collision dynamics of the solitary waves as solutions of linearly Coupled Nonlinear Schr\"odinger Equations for various initial phases. The initial conditions are superposition of two one-soliton solutions with general polarization. The quasi-particle behavior of propagating and interacting solutions in conditions of rotational polarization is examined. We find that the total mass, pseudomomentum and energy are conserved while the local masses, individual and total polarization depend strongly on the linear coupling and the initial phase difference. We also find out that the polarization angle of the quasi-particles can change independently of the interaction.
\end{abstract}

\maketitle

\section{Introduction}

The system of Coupled Nonlinear Schr{\"o}dinger Equations (CNLSE) is an attractive because of several reasons. From among them: CNLSE are a soliton supporting system, the system is a proper ``a spring board'' for investigate and track the quasi-particle (QP) behavior of the solitary waves; it possesses very rich phenomenology which can be extended adding new terms in the equations; the vector structure of the system and their solutions admits to study and treat the polarization vector \cite{lakoba}, \cite{yang}, \cite{hempelmann}. The initial polarization of  CNLSE and its evolution during the QP interaction is a very  important element. There are numerous experimental observations (see \cite{collings} and cited there literature) of the effects connected with this quantity. For the full fledged CNLSE with general initial polarization, no analytic solution is available. It can be obtained, e.g., via adequately devised numerical scheme. An implicit scheme of Crank-Nicolson type was first proposed for the single NLSE in \cite{TahaAblovitz}. In order to adapt it to CNSLE, internal iterations were applied in \cite{ChriDostMau}. The latter yielded two very important advantages: fully implicit scheme and implementation of the conservation laws on difference level within the round-off error of the calculations. The above scheme was extended to complex arithmetic  in \cite{TodoChri06} where the computer code for Gaussian elimination with pivoting of \cite{Chri_Reading} was generalized for complex-valued  multi-diagonal band algebraic systems. Having in mind the importance of initial polarization, we investigated in \cite{TodoChri08} the collision dynamics for circularly polarized solitons based on $sech$-functions and found out that there exist infinite number of Manakov two-soliton solutions preceded by polarization discontinuity (shock) on the place of interaction. In order to extend the range of investigations in the case of general polarization we established an auxiliary conjugate system of nonlinear ordinary differential equations \cite{bgsiam09} in order to generate numerically initial elliptically polarized soliton solutions and applied them for numerous head-on \cite{matcom10} and takeover collisions \cite{TodoAIMS10}. To enrich the concept of the polarization we involved into considerations linear coupling of $sech$-like solitons which is a generator of their rotational polarization (breathing) and/or gain/self-dissipation \cite{aip1340}. In this work we aim to conduct series of simulations and to track the particle-like behavior of the interacted localized non-$sech$ waves with arbitrary (elliptic) initial polarization.

\section{Problem formulation}
CNLSE is system of nonlinearly coupled Schr\"odinger equations (also called the Gross-Pitaevskii or Manakov-type system):
\begin{subequations}\label{strong}
\begin{align}
{\rm i} \psi_t &= \beta\psi_{xx} + \bigl[\alpha_1|\psi|^2 +
(\alpha_1+2\alpha_2)|\phi|^2\bigr]\psi ,
\label{strongfirst} \\
{\rm i} \phi_t &= \beta\phi_{xx} + \bigl[\alpha_1|\phi|^2 +
(\alpha_1+2\alpha_2)|\psi|^2\bigr]\phi,\label{strongsecond}
\end{align}
\end{subequations}
where $\beta$ is the dispersion coefficient, $\alpha_1$ describes the self-focusing of a signal for pulses in birefringent
media, $\alpha_2$ (called cross-modulation parameter) governs the nonlinear coupling between the equations. When $\alpha_2=0$, no nonlinear coupling is present despite the fact that ``cross-terms'' proportional to $\alpha_1$ appear in the equations and the both equations admit identical solution equal to the solution of single  NLSE with nonlinearity coefficient $\alpha=2\alpha_1$. Now we turn our considerations to a modified CNLSE where the nonlinear coupling is trivial but is replaced by so-called linear coupling. Nevertheless the system remains nonlinear
\begin{subequations}\label{eq:CNLSE}
\begin{align}
{\rm i} \psi_t &= \beta\psi_{xx} + \alpha_1\bigl(|\psi|^2 + |\phi|^2\bigr)\psi-\Gamma \phi, \label{eq:CNLSE_psig}\\
{\rm i} \phi_t &= \beta\phi_{xx} + \alpha_1\bigl(|\phi|^2 + \psi|^2\bigr)\phi-\Gamma \psi. \label{eq:CNLSE_phig}
\end{align}\end{subequations}
The magnitude of linear coupling is governed by the complex-valued quantity $\Gamma$. $\Re\Gamma$ governs the oscillations between states termed as breathing solitons and actually seems their frequency of oscillation, while $\Im\Gamma$ describes the gain/self-dissipation behavior of soliton solutions. As it is shown in \cite{Optical_IMACS} the breathing is present free of the interaction. The relation between the Manakov-type system (\ref{strong}) (when $\alpha_2=0$) and the Linearly CNLSE (\ref{eq:CNLSE}) is given by the substitution (for details see \cite{Optical_IMACS} and \cite{aip1340})
\begin{equation}
\psi = \Psi \cos(\Gamma t) + {\rm i} \Phi \sin(\Gamma t),\quad \phi = \Phi \cos(\Gamma t) + {\rm i} \Psi \sin(\Gamma t),\label{eq:lCNLSE}
\end{equation}
where $\Phi$ and $\Psi$ are soliton solutions  of \eqref{strong}. These properties and references \eqref{eq:lCNLSE} of the solutions allow and argue us to solve (\ref{eq:CNLSE}) instead (\ref{strong}) in order to get more information about the breathing behavior of the soliton solutions. Let us remind that in \cite{aip1340} we considered the problem in question but in the particular case when $\Phi$ and $\Psi$ were assumed to be certainly \textit{sech}-solutions of \eqref{strong}.

In this paper the solutions are pulses whose modulation amplitude is of general form (non-\emph{sech}) and their polarization\footnote{The latter determines one more kind of polarization called rotational polarization.} rotates with time. This determines the choice of initial conditions for numerical investigation of temporal evolution of interacting solitons. As usual we concern ourselves with the soliton solutions which are localized envelopes on a propagating carrier wave.
We assume that for each of the functions $\phi,\psi$ the initial condition has the general type
\begin{equation}
\vec\chi (x,t;X,c,n_\chi) = A_{\chi}(x+X-ct)\exp{\left\{{\rm i}\left[n_{\chi} t- \tfrac{1}{2}c(x-X-ct)+\delta_{\chi}\right]\right\}}\label{eq:incond}
\end{equation}
where $\vec\chi$ stands for $(\psi,\phi)$, $c$ is the phase speed of the envelope, $X$ is the initial position of the center of the soliton; $\vec{n} (n_{\psi}$, $n_{\phi})$ is the vector of carrier frequencies of the components; $\vec{\delta}(\delta_{\psi},\delta_{\phi})$ is phase vector of the two components.
Note that the phase speed $c$ must be the same for the two components $\psi$ and $\phi$. Otherwise they will propagate with different phase speeds and after some time the two components will be in two different positions in space, and will no longer form a single structure. For the envelopes $(A_\psi, A_\phi)$, $\theta\equiv \arctan ({\max |\phi|}/{\max |\psi|})$  is a polarization angle. We are interested in solutions with initially nontrivial carrier frequencies ($n_{\psi}, n_{\phi}) \neq (0,0)$, i.e., with elliptic or circular polarizations. Let us emphasize that these qualities are intrinsic to non-\emph{sech}-like solutions.
Having in mind (\ref{eq:incond}) we implement straightforward manipulations yielding
\begin{subequations}
\begin{align}
&\left[A_{\psi}''+\big(n_{\psi}+\tfrac{1}{4}c^2\big)A_{\psi}+\alpha_1 \left(A_{\psi}^2+ A_{\phi}^2\right)A_{\psi}\right]\exp{\left\{{\rm i}\left[n_{\psi} t- \tfrac{1}{2}c(x-X-ct)+\delta_{\psi}\right]\right\}}\cos(\Gamma t)\notag \\+ &{\rm i} \left[A_{\phi}''+\big(n_{\phi}+\tfrac{1}{4}c^2\big)A_{\phi}+\alpha_1 \left(A_{\phi}^2+ A_{\psi}^2\right)A_{\phi}\right]\exp{\left\{{\rm i}\left[n_{\phi} t- \tfrac{1}{2}c(x-X-ct)+\delta_{\phi}\right]\right\}}\sin(\Gamma t)=0,\\
&\left[A_{\phi}''+\big(n_{\phi}+\tfrac{1}{4}c^2\big)A_{\phi}+\alpha_1 \left(A_{\phi}^2+ A_{\psi}^2\right)A_{\phi}\right]\exp{\left\{{\rm i}\left[n_{\phi} t- \tfrac{1}{2}c(x-X-ct)+\delta_{\phi}\right]\right\}}\cos(\Gamma t)\notag \\+ &{\rm i} \left[A_{\psi}''+\big(n_{\psi}+\tfrac{1}{4}c^2\big)A_{\psi}+\alpha_1 \left(A_{\psi}^2+ A_{\phi}^2\right)A_{\psi}\right]\exp{\left\{{\rm i}\left[n_{\psi} t- \tfrac{1}{2}c(x-X-ct)+\delta_{\psi}\right]\right\}}\sin(\Gamma t)=0.
\end{align}
\end{subequations}
The above two equations are true if and only if the two components of the envelope are governed by the following coupled system of differential equations:
\begin{subequations}\label{eq:conjsys}
\begin{align}
A_{\psi}''+\big(n_{\psi}+\tfrac{1}{4}c^2\big)A_{\psi}+\alpha_1 \left(A_{\psi}^2+ A_{\phi}^2\right)A_{\psi}=0, \label{eq:conjsys_psig}  \\[1mm]
A_{\phi}''+\big(n_{\phi}+\tfrac{1}{4}c^2\big)A_{\phi}+\alpha_1\left( A_{\phi}^2+A_{\psi}^2\right)A_{\phi}=0.  \label{eq:conjsys_phig}
\end{align}
\end{subequations}
  \begin{figure}[ht!]
\centerline{\includegraphics[width=0.75\textwidth]{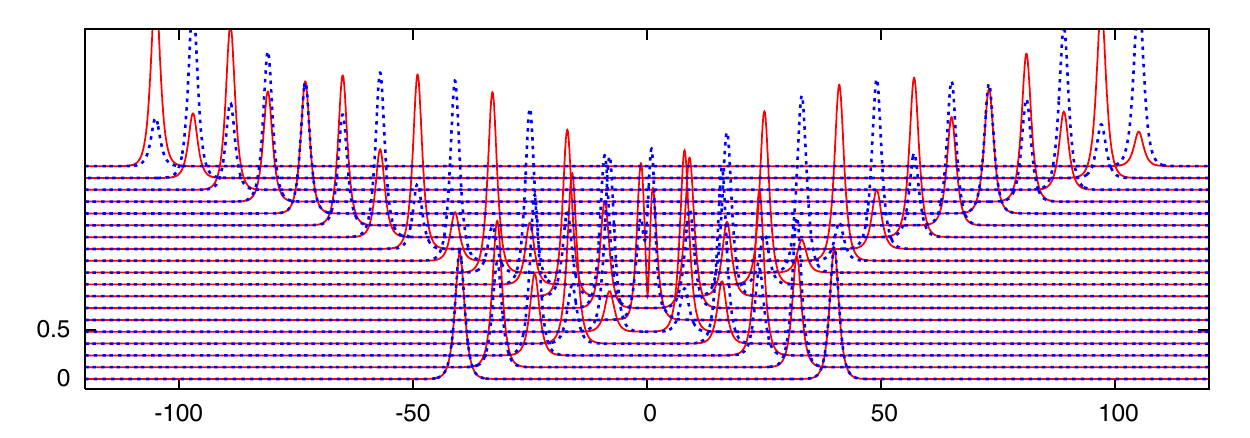}}
\vspace{-0.1in}
\centerline{\small (a) $\deljump=0^\circ$}
\centerline{\includegraphics[width=0.75\textwidth]{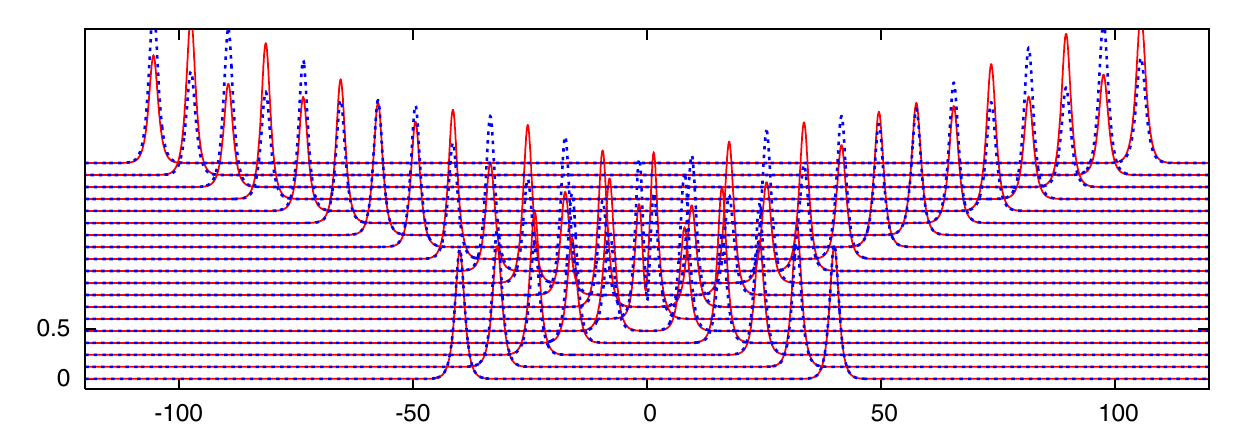}}
\vspace{-0.1in}
\centerline{\small (b) $\deljump=90^\circ$}
\centerline{\includegraphics[width=0.75\textwidth]{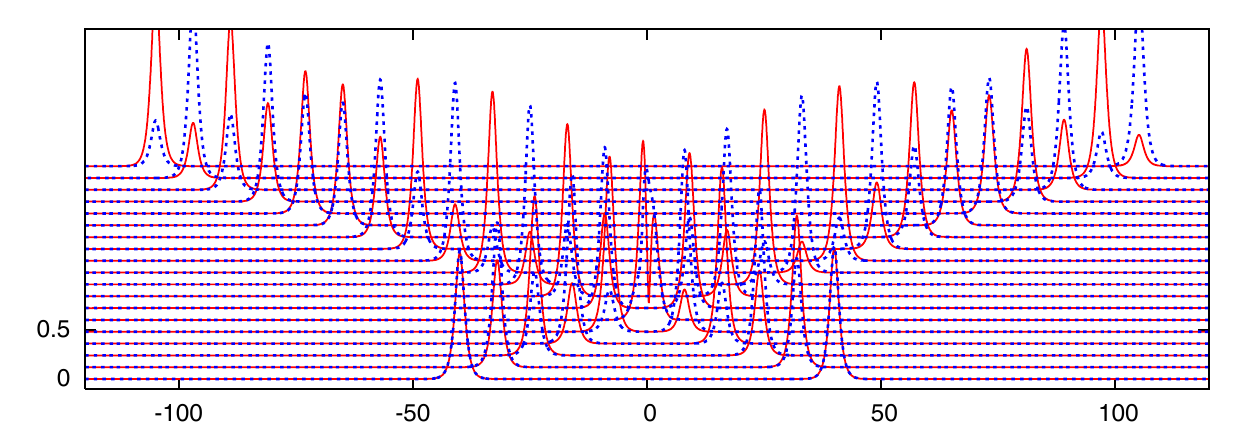}}
\vspace{-0.1in}
\centerline{\small (c) $\deljump=180^\circ$}
\caption{Head-on collision of QPs with initial circular polarization for $n_{l \psi}=n_{r \psi}=n_{l \phi}=n_{r \phi}=-1.5$, $c_l=-c_r=1$, $\alpha_1=0.75$, $\Gamma=0.175$, $X_l=-40$, $X_r=40$ and different initial phase differences.}
\label{fig:clashcird180}
\end{figure}
We find out that this system is explicitly equivalent to the bifurcation conjugate system\footnote{The system admits bifurcation solutions since the trivial solution obviously is always present.} corresponding to Manakov case with $\alpha_2=0$ (see \cite{matcom10} and \cite{bgsiam09}). Here we aim to trace the strong linear coupling on the interaction in conditions of general initial polarization. After that we construct the initial conditions for $\vec\chi$ from (\ref{eq:incond}) as a superposition of two two-component solitary waves located at positions $X_l$ and $X_r$ and propagating with phase speeds $c_l$ and $c_r$
\begin{equation}
\vec\chi(x,0)=\vec\chi_l(x,0,c_l,n_l,X_l,\vec\delta_l)+\vec\chi_r(x,0,c_r,n_r,X_r,\vec\delta_r).
\label{eq:incnd_superpos}
\end{equation}
\begin{figure}[ht!]
\centerline{\includegraphics[width=0.75\textwidth]{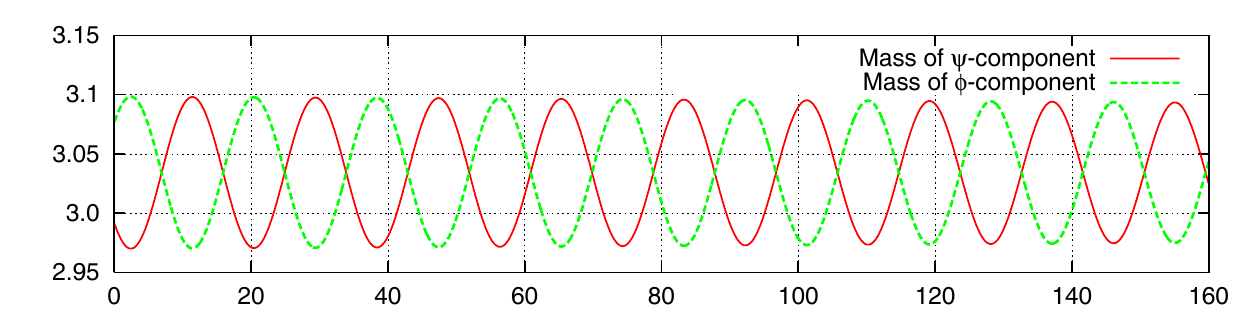}}
\centerline{\small (a) $\deljump=0^\circ$}
\centerline{\includegraphics[width=0.75\textwidth]{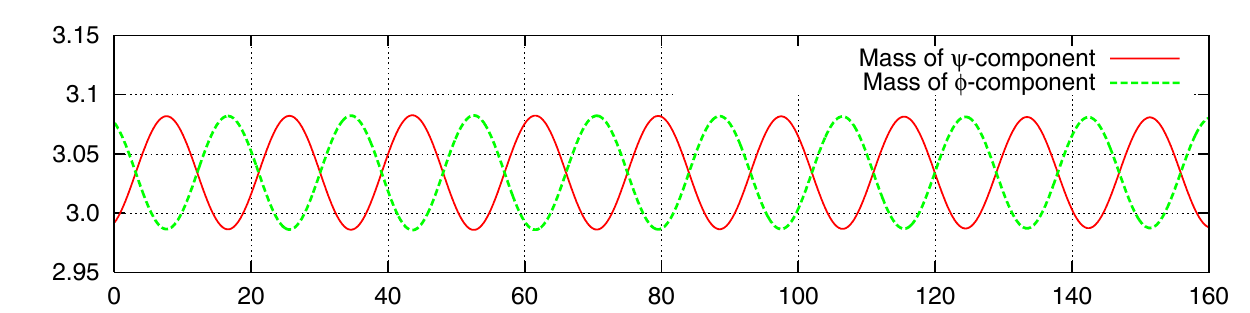}}
\centerline{\small (b) $\deljump=90^\circ$}
\centerline{\includegraphics[width=0.75\textwidth]{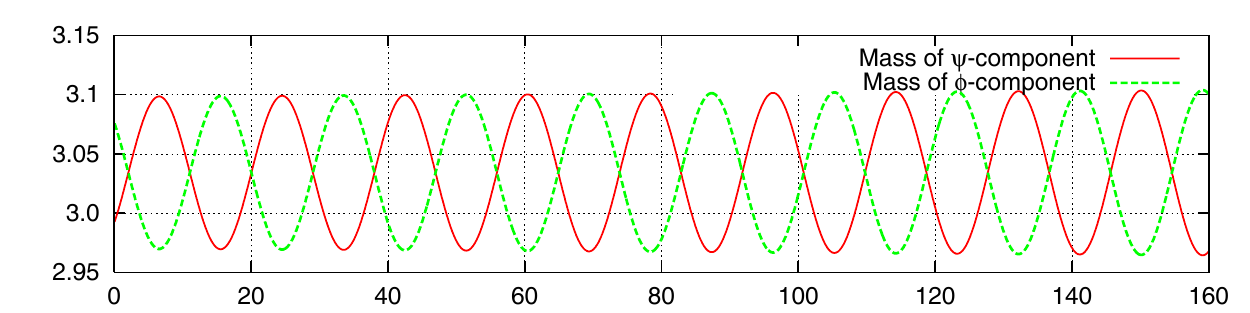}}
\centerline{\small (c) $\deljump=180^\circ$}
\caption{ Mass dynamics corresponding to Figure~\ref{fig:clashcird180}.}
\label{fig:massmomcird90}
\end{figure}
\begin{figure}[ht!]
\centerline{\includegraphics[width=0.75\textwidth]{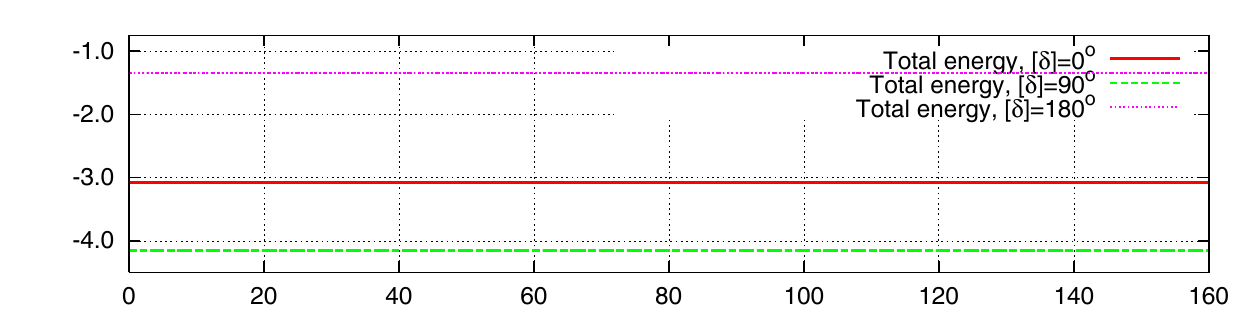}}
\vspace{-0.1in}
\caption{Influence of the initial phase difference on the energy dynamics corresponding to Figure~\ref{fig:clashcird180}.}
\label{fig:energycird0}
\end{figure}
\begin{figure}[ht!]
\centerline{\includegraphics[width=0.75\textwidth]{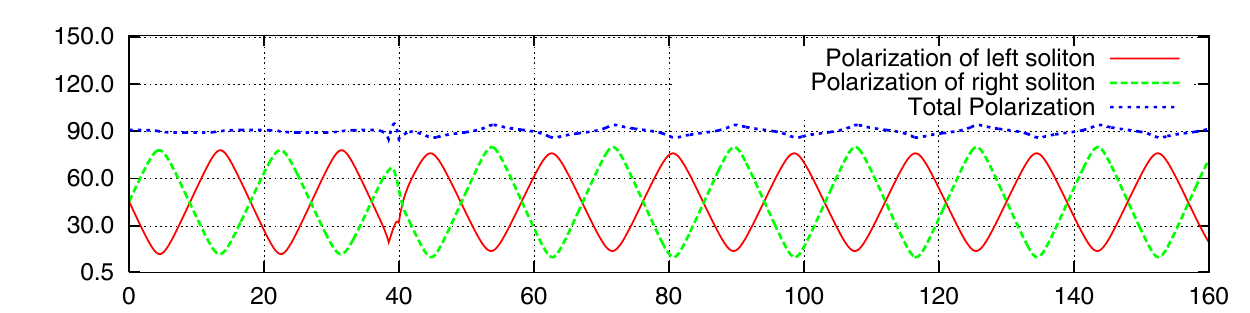}}
\centerline{\small (a) $\deljump=0^\circ$}
\centerline{\includegraphics[width=0.75\textwidth]{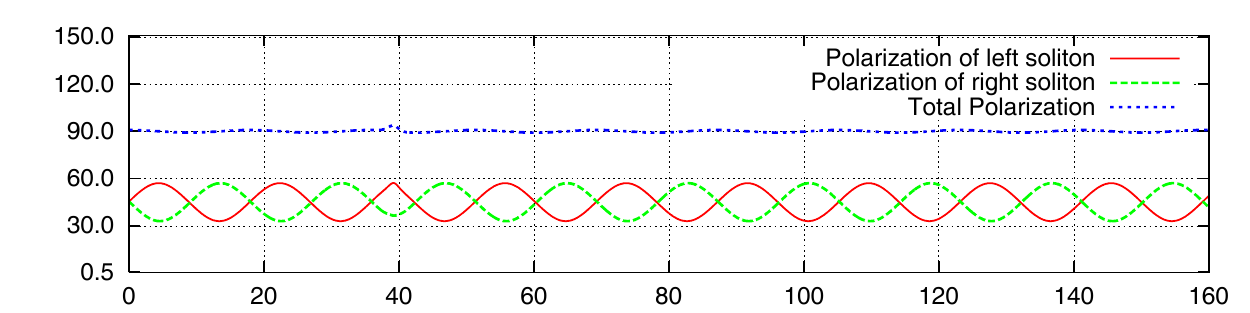}}
\centerline{\small (b) $\deljump=90^\circ$}
\centerline{\includegraphics[width=0.75\textwidth]{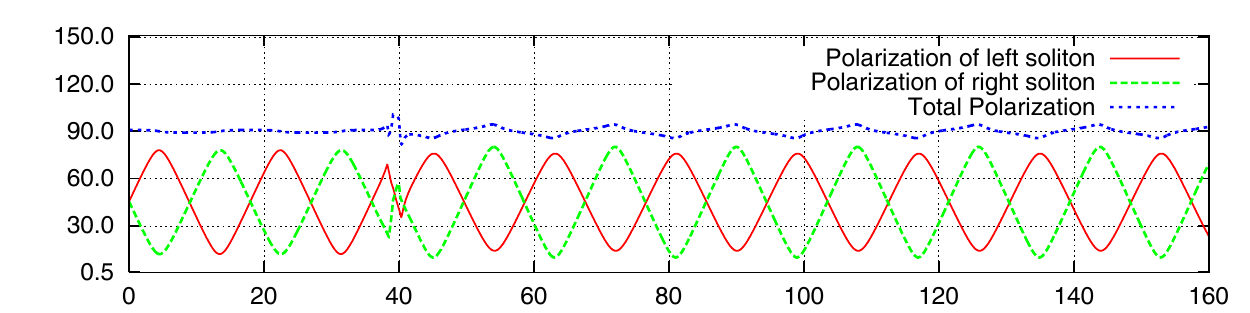}}
\centerline{\small (c) $\deljump=180^\circ$}
\caption{ Individual and total polarization dynamics in dependence on initial phase difference corresponding to Figure~\ref{fig:clashcird180}.}
\label{fig:polarcird180}
\end{figure}
The initial distance $|X_l-X_r|$ is supposed to be large enough in order the solutions not to overlay in the time moment $t=0$.
The system (\ref{eq:conjsys}) is solved numerically and their solutions are  readily used as an approximate initial condition for the unsteady computations.
Our aim is to understand better the influence of the
initial polarization and the initial phase difference on the particle-like behavior of
the localized waves (QPs). They survive the collision with other QPs (or some other
kind of interactions) without losing its identity. The initial difference  in phases can have a profound influence on the polarizations of the solitons after the interaction as well as on the magnitude of the full energy and the amplitude of breathing.

Before turning to the numerical investigation we mention here that the system Eq.~\eqref{eq:CNLSE} possesses three conservation laws when asymptotic boundary conditions are imposed, namely when $\psi,\phi \rightarrow 0$ for $x \rightarrow \pm\infty$.  Following \cite{ChriDostMau,Optical_IMACS} we define  ``mass'', $M$, (pseudo)momentum, $P$, and energy, $E$  as follows
\begin{equation}
M \stackrel{\mathrm{def}}{=} \frac{1}{2\beta}\int_{-L_1}^{L_2}\!\!\!\left( |\psi|^2 + |\phi|^2 \right)\! {\rm d}x,  \quad
P \stackrel{\mathrm{def}}{=}  -\int_{-L_1}^{L_2} \!\!\!{\cal {I}}(\psi \bar\psi_x + \phi \bar \phi_x) \mathrm{d}x, \quad
E \stackrel{\mathrm{def}}{=} \int_{-L_1}^{L_2} \!\!\!\mathcal{H} \mathrm{d}x,
\label{eq:5}
\end{equation}
where\vspace{-0.2in}
\begin{equation*}
\mathcal{H} \stackrel{\mathrm{def}}{=}\beta \left(
|\psi_x|^2 + |\phi_x|^2\right) - \frac{\alpha_1}{2}(|\psi|^2 +
|\phi|^2)^2 -2\Gamma[\Re({\bar\psi}{\bar\phi})]
\end{equation*}
is the Hamiltonian density of the system. Here $-L_1$ and $L_2$ are the left end and the right end of the interval under consideration.  The conservation laws read
\begin{equation}
\frac{d M}{d t} = 0,\qquad \frac{d P}{d t} = 0,\qquad \frac{d E}{d t } = 0,
\label{eq:conserv_laws}
\end{equation}
which means that for asymptotic boundary conditions CNLSE admit at most 3 conservation laws, i.e.,  the system \eqref{eq:CNLSE} is non-fully integrable.

\section{Numerical method}

 In order to be able to obtain reliable results for the time evolution of the solution, one needs to devise a difference scheme that represent faithfully the above conservation laws. Such a scheme was proposed in \cite{ChriDostMau}, and applied in \cite{Optical_IMACS}. This scheme was based on a fast Gaussian elimination solver for multi-diagonal systems \cite{Chri_Reading}.  Consequently,  the above mentioned scheme was implemented for complex arithmetic  in \cite{TodoChri06}, where  a complex arithmetic algorithms was developed to generalize the one from \cite{Chri_Reading}. The complex-arithmetic algorithm is four times faster, and we will use it also in the present paper. Thus,  for solving  Eqs.~\eqref{eq:CNLSE} with the initial conditions \eqref{eq:incond}, \eqref{eq:incnd_superpos} numerically, we use an implicit conservative scheme in complex arithmetic:
 \begin{subequations}\label{eq:non_lin_sche}
  \begin{multline}
 {\rm i} \frac{\psi^{n+1}_i - \psi^n_i}{\Delta\tau} = \frac{\beta}{2h^2}
\left(\psi^{n+1}_{i-1} - 2\psi^{n+1}_{i} + \psi^{n+1}_{i+1} +
\psi^{n}_{i-1} - 2\psi^{n}_{i} + \psi^{n}_{i+1}\right) \\
 + \frac{\alpha_1}{4}\left(\psi^{n+1}_i + \psi^n_i\right) \Bigl[\left(|\psi^{n+1}_{i}|^2 + |\psi^{n}_{i}|^2 \right)+ \left(|\phi^{n+1}_{i}|^2 + |\phi^{n}_{i}|^2 \right)\Bigr] - \frac{\Gamma}{2} \left(\phi_i^{n+1}+\phi_i^n\right),\label{eq:non_lin_sche_psi}
 \end{multline}\vspace{-0.3in}
 \begin{multline}
{\rm i} \frac{\phi^{n+1}_i - \phi^n_i}{\Delta\tau} = \frac{\beta}{2h^2} \left(\phi^{n+1}_{i-1} - 2\phi^{n+1}_{i} + \phi^{n+1}_{i+1} + \phi^{n}_{i-1} - 2\phi^{n}_{i} + \phi^{n}_{i+1}\right) \\
 + \frac{\alpha_1}{4} \left(\phi^{n+1}_i + \phi^n_i\right) \Bigl[\left(|\phi^{n+1}_{i}|^2 + |\phi^{n}_{i}|^2 \right)+\left(|\psi^{n+1}_{i}|^2 + |\psi^{n}_{i}|^2 \right)\Bigr]- \frac{\Gamma}{2} \left(\psi_i^{n+1}+\psi_i^n\right),\label{eq:non_lin_sche_phi}
 \end{multline}
\end{subequations}
\begin{figure}[ht!]
\centerline{\includegraphics[width=0.68\textwidth]{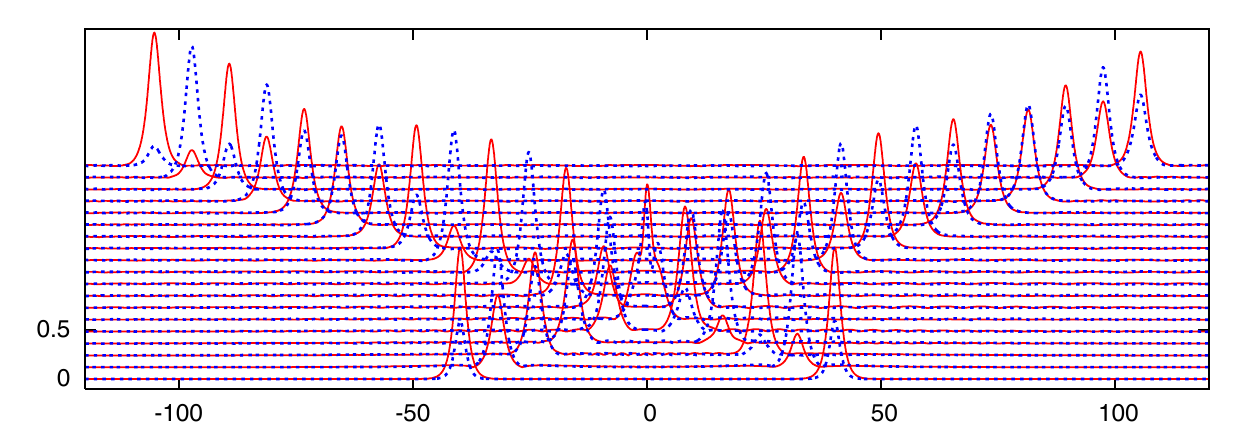}}
\centerline{\small (a) $\deljump=0^\circ$}
\centerline{\includegraphics[width=0.68\textwidth]{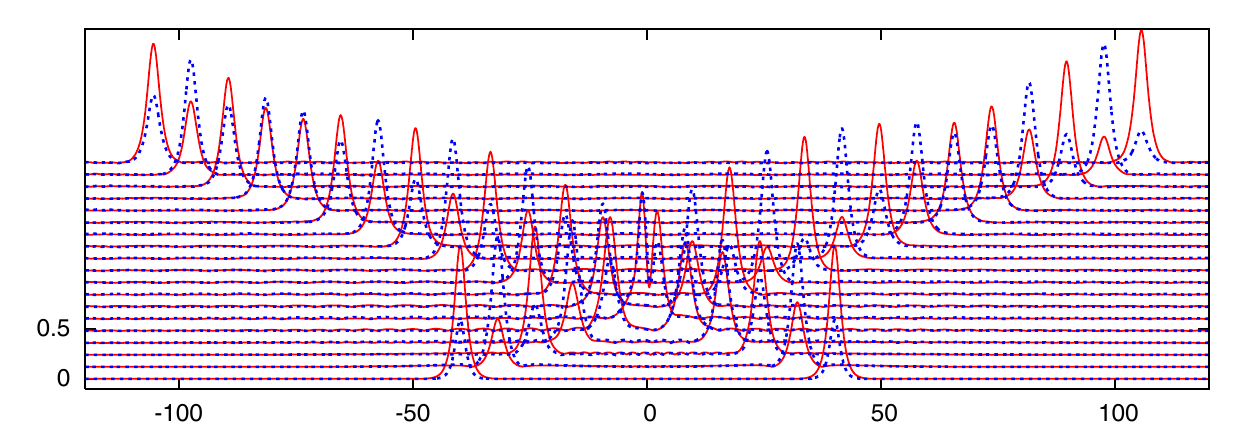}}
\centerline{\small (b) $\deljump=90^\circ$}
\centerline{\includegraphics[width=0.68\textwidth]{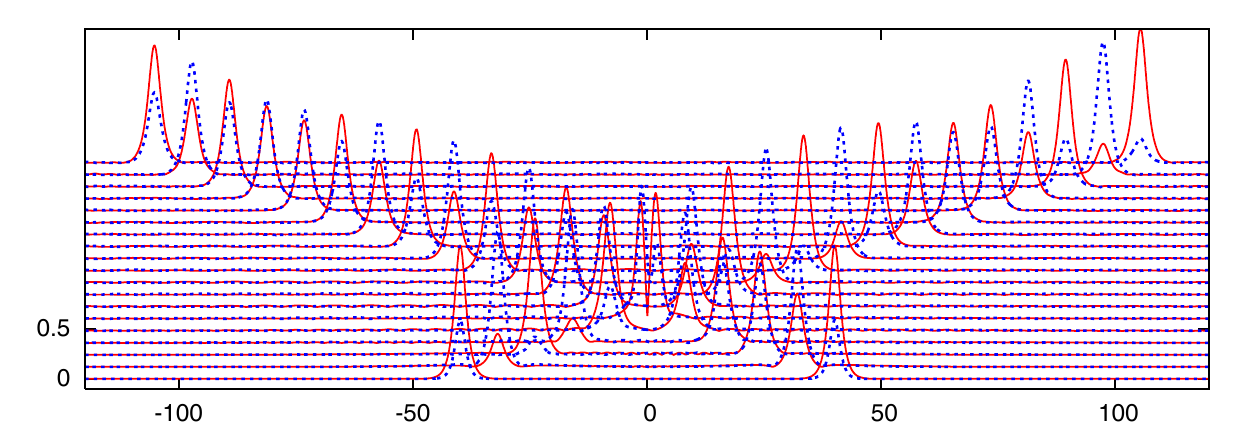}}
\centerline{\small (c) $\deljump=135^\circ$}
\centerline{\includegraphics[width=0.68\textwidth]{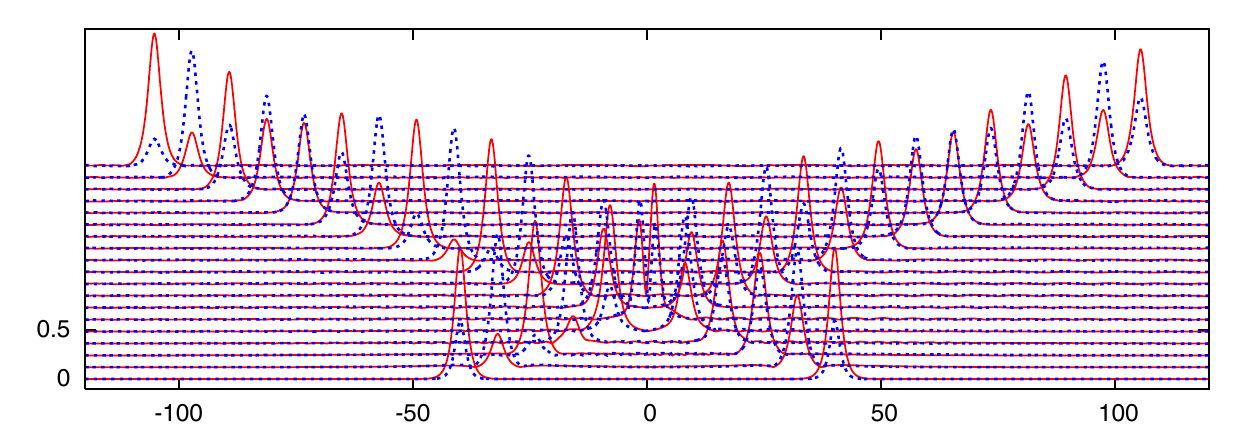}}
\centerline{\small (d) $\deljump=180^\circ$}
\caption{Head-on collision of QPs with initial elliptic polarization for $\theta=23^\circ 44'$, $n_{l \psi}=n_{r \psi}=-1.1$, $n_{l \phi}=n_{r \phi}=-1.5$, $c_l=-c_r=1$, $\alpha_1=0.75$, $\Gamma=0.175$, $X_l=-40$, $X_r=40$ and different initial phase differences.}
\label{fig:clasheld180}
\end{figure}
\begin{figure}[ht!]
\centerline{\includegraphics[width=0.75\textwidth]{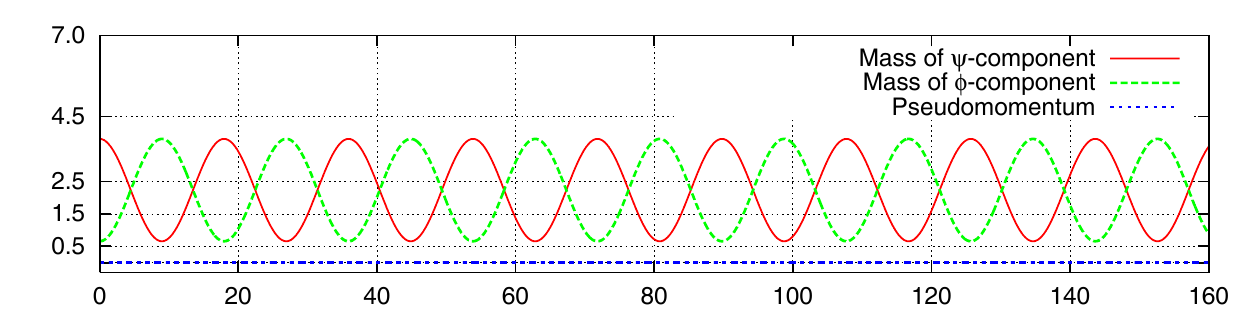}}
\caption{Mass and pseudomomentum dynamics for $\deljump=0^\circ$ corresponding to Figure~\ref{fig:clasheld180}.}
\end{figure}
\begin{figure}[h!]
\centerline{\includegraphics[width=0.75\textwidth]{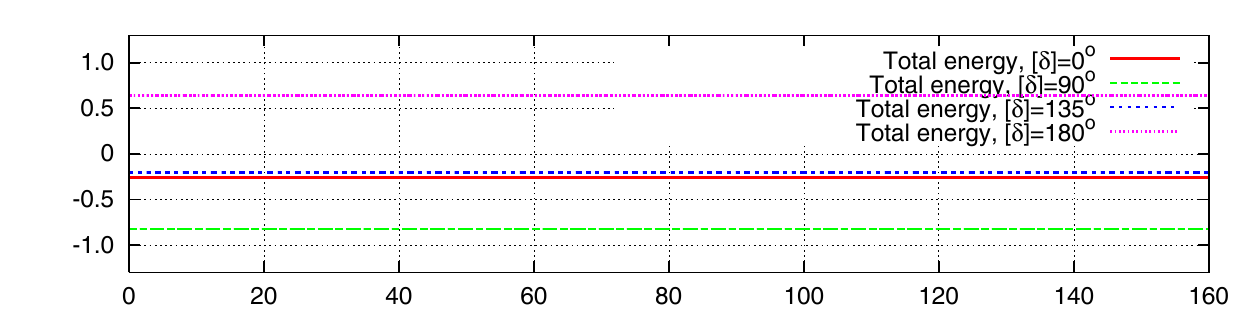}}
\vspace{-0.1in}
\caption{ Influence of the initial phase difference on the total energy dynamics corresponding to Figure~\ref{fig:clasheld180}.}
\label{fig:ell_allenergiies}
\end{figure}
\begin{figure}[th!]
\centerline{\includegraphics[width=0.75\textwidth]{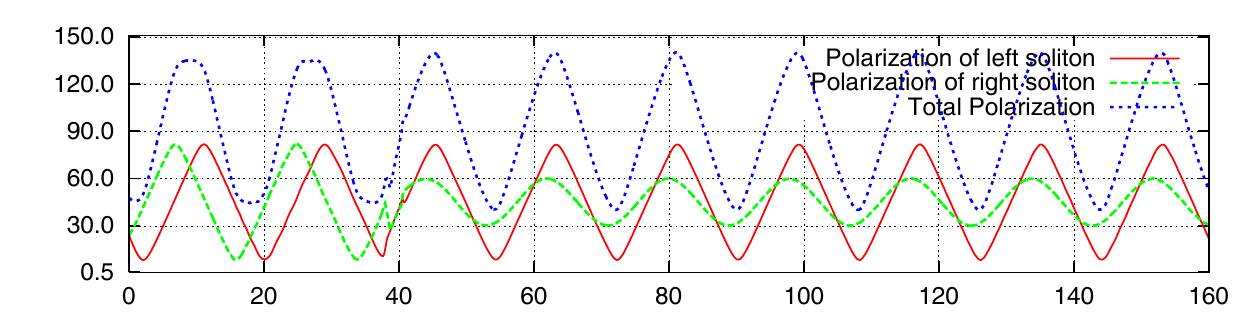}}
\centerline{\small (a) $\deljump=0^\circ$}
\centerline{\includegraphics[width=0.75\textwidth]{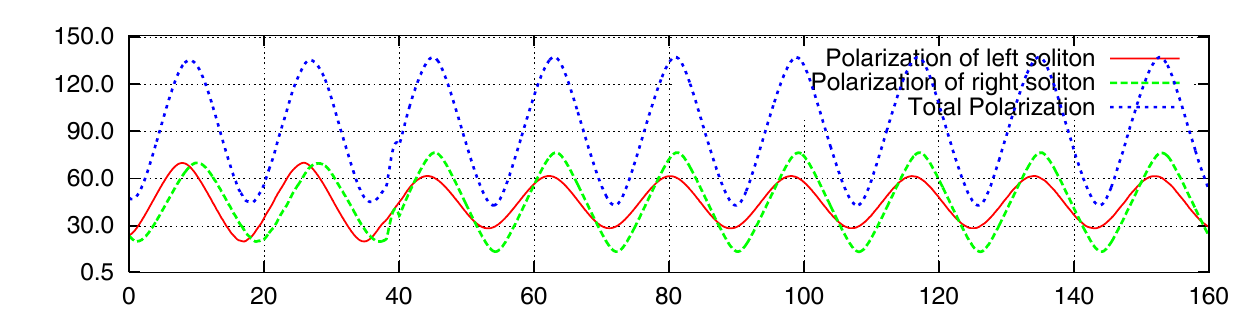}}
\centerline{\small (b) $\deljump=90^\circ$}
\centerline{\includegraphics[width=0.75\textwidth]{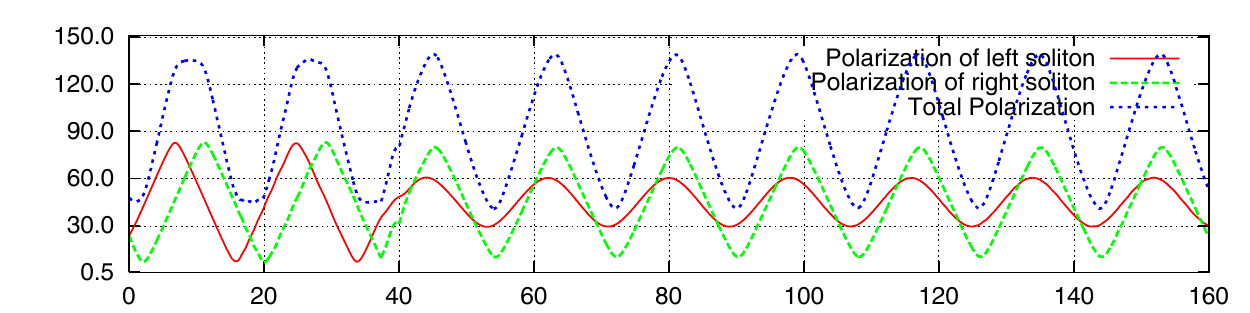}}
\centerline{\small (c) $\deljump=135^\circ$}
\centerline{\includegraphics[width=0.75\textwidth]{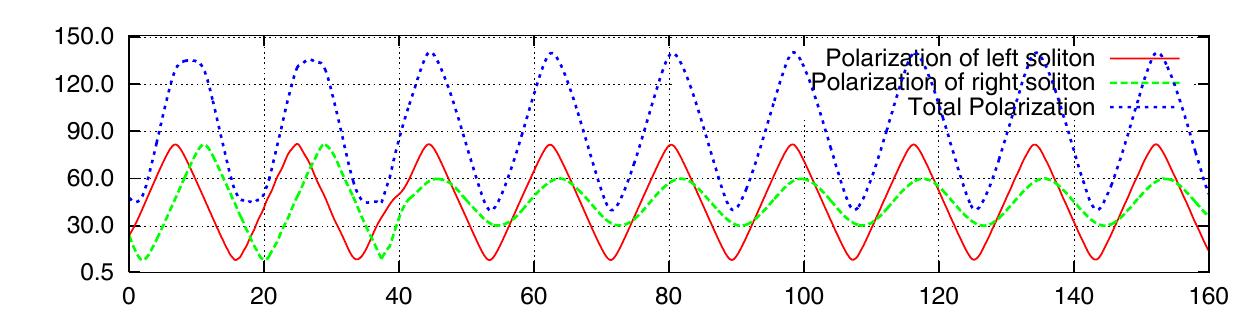}}
\centerline{\small (d) $\deljump=180^\circ$}
\caption{Individual and total polarization dynamics in dependence on initial phase difference corresponding to Figure~\ref{fig:clasheld180}.}
\label{fig:polareld180}
\end{figure}
\noindent on mesh $(x_i,\tau^n)$ with $x_i=-L_1+i h$, $h = {(L_2+L_1)}/{m}$, $i=1,...,m$ and $\tau^n = n \Delta \tau$, $n=0,1,2...$
It is not only convergent (consistent and stable), but also conserves mass, pseudomomentum, and energy, i.e., there exist discrete analogs $M^n$, $P^n$ and $E^n$,  for \eqref{eq:5}, which arise from the scheme (for details see \cite{ChriDostMau,Optical_IMACS,TodoChri06}).
\begin{align}
M^n &= \sum_{i=2}^{N-1} \left(|\psi^n_i|^2 + |\phi^n_i|^2 \right)= {\rm const}, \nonumber \\
P^n &= -\sum_{i=2}^{N-1}\frac{1}{h} \Im \left[\psi^n_i \left(\bar\psi_{i+1}^n-\bar\psi_i^n\right) + \phi_i^n \left(\bar\phi_{i+1}^n - \bar\phi_i^n\right)\right]= {\rm const}, \nonumber \\
E^n&=\sum_{i=2}^{N-1}\frac{-\beta}{2h^2}\left(|{\psi}^{\;n}_{i+1}-{\psi}^{\;n}_i|^2+ |{\phi}^{\;n}_{i+1}-
{\phi}^{\;n}_i|^2\right)+ \frac{\alpha_1}{4}\left(|\psi^{n}_i|^2+|\phi^{n}_i|^2\right)^2+ \Gamma {\Re}[\bar \phi_i^n \psi_i^n]
= {\rm const}, \nonumber
\end{align}
for $n\geq 0$. These values are kept constant by the scheme during the time stepping. The above scheme is of  Crank-Nicolson type for the linear terms and we employ internal iteration to achieve implicit approximation of the nonlinear terms, i.e., we use its linearized implementation \cite{ChriDostMau}. In this way the order of approximation of Eqs.~\eqref{eq:non_lin_sche} is $O(\Delta\tau^2+h^2)$.

The above nonlinear scheme is implemented via internal iterations as follows:
\begin{subequations}
\begin{align}\label{eq:psi_scheme1}
{\rm i} \frac{\psi^{n+1,k+1}_i - \psi^n_i}{\Delta\tau} =
&\frac{\beta}{2h^2} \Big(\psi^{n+1,k+1}_{i-1} -
2\psi^{n+1,k+1}_{i} + \psi^{n+1,k+1}_{i+1}
+ \psi^{n}_{i-1} -  2\psi^{n}_{i} + \psi^{n}_{i+1}\Big) \nonumber \\ &\hspace{-1.2in}+ \frac{\alpha_1}{4}\left(\psi^{n+1,k}_i+ \psi^n_i\right)  \Big[\big(|\psi^{n+1,k+1}_{i}||\psi^{n+1,k}_{i}|  + |\psi^{n}_{i}|^2\big)
+ \big(|\phi^{n+1,k+1}_{i}| |\phi^{n+1,k}_{i}| + |\phi^{n}_{i}|^2 \big)\Big]\nonumber\\ &- \frac{\Gamma}{2}(\psi_i^n+\psi_i^{n+1,k}),\\
{\rm i} \frac{\phi^{n+1,k+1}_i - \phi^n_i}{\Delta\tau} = &\frac{\beta}{2h^2} \big(\phi^{n+1,k+1}_{i-1} - 2\phi^{n+1,k+1}_{i} + \phi^{n+1,k+1}_{i+1} + \phi^{n}_{i-1} -2\phi^{n}_{i} + \phi^{n}_{i+1}\big) \nonumber \\ &\hspace{-1.2in}+ \frac{\alpha_1}{4} \left(\phi^{n+1,k}_i + \phi^n_i\right)
\Big[\big(|\phi^{n+1,k+1}_{i}||\phi^{n+1,k}_{i}| + |\phi^{n}_{i}|^2 \big)
+\big(|\psi^{n+1,k+1}_{i}||\psi^{n+1,k}_{i}| + |\psi^{n}_{i}|^2 \big)\Big]\nonumber\\ &- \frac{\Gamma}{2}(\phi_i^n+\phi_i^{n+1,k}). \label{eq:phi_scheme2}
\end{align}
\end{subequations}

The iterations are repeated (stepping up the index `$k$') until convergence is reached. If the internal iterations are convergent, the scheme in full steps Eqs.~\eqref{eq:non_lin_sche} is absolutely stable due to its conservative nature.  For not very large time steps the iteration requires less than half a dozen loops to converge.  This is a small price to pay having in mind the inextricably coupled five-diagonal complex structure of the matrix.

The time-stepping scheme needs initial conditions of the type of Eq.~\eqref{eq:incond}.  In the general case of elliptic initial polarization we discretize the system Eqs.~\eqref{eq:conjsys} in the same manner as the original evolutionary system, and use Newton's method and hermitian  splines to get approximately the solution \cite{bgsiam09}. The result are used as initial conditions for the problem under consideration.

The above presented scheme and algorithm have been verified  for different girds and time increments and the approximation has been confirmed.

%

\section{Results and discussion}
\subsection{Circular polarization. Head-on collision}
This is a special elliptic polarization with polarization angle $\theta=45^\circ$. The initial configuration is generated from the auxiliary bifurcation system (\ref{eq:conjsys}). Because the parametric space of the problem is too big to be explored in full, we choose $n_{l \psi}=n_{r \psi}=n_{l \phi}=n_{r \phi}=-1.5$, $c_l=-c_r=1$, $\alpha_1=0.75$, $\Gamma=0.175$ and focus on the effects of $\deljump \equiv \delta_r-\delta_l$. We observe that the breathing does not interfere with the soliton collision since they are distinguished from previously obtained breathers \cite{1}, \cite{2}. One sees the ``breathing" of the pulses even without any interaction. This is the manifestation of the rotation of the polarization. Thus $\Gamma$ is responsible for the exchange of wave mass between the modes and we call it ``cross" dispersion of the signals (see Figure~\ref{fig:massmomcird90}). \begin{figure}[ht!]
\centerline{\includegraphics[width=0.75\textwidth]{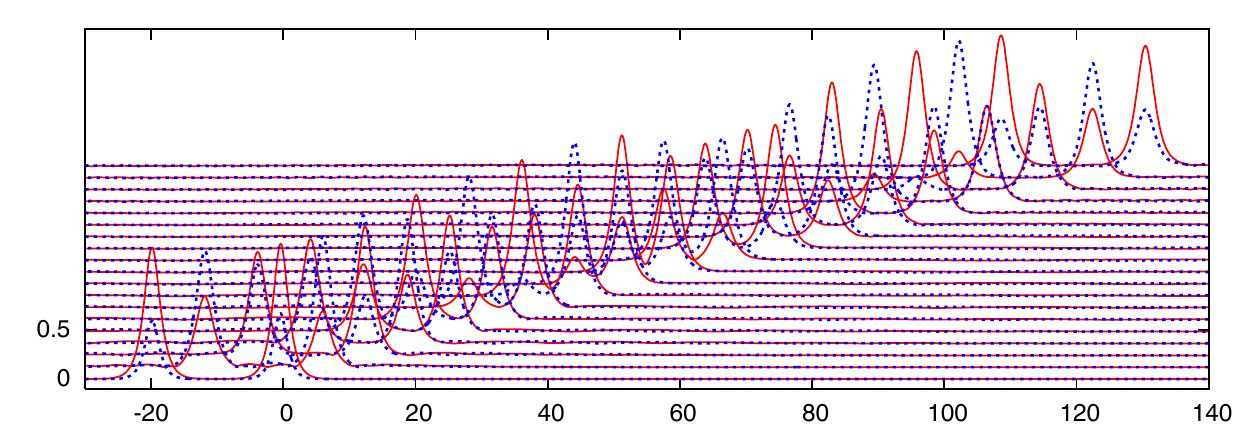}}
\centerline{\small (a) $\deljump=0^\circ$}
\centerline{\includegraphics[width=0.75\textwidth]{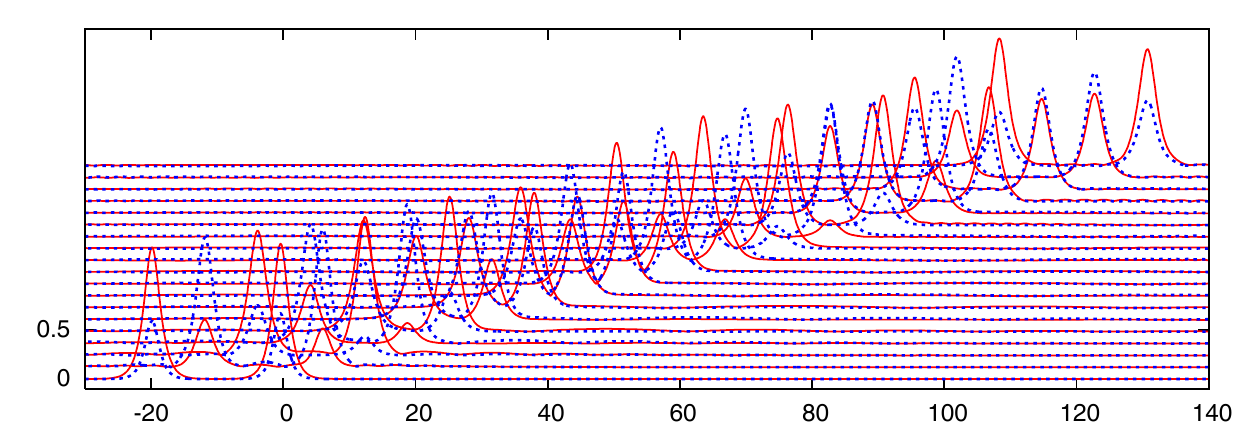}}
\centerline{\small (b) $\deljump=90^\circ$}
\centerline{\includegraphics[width=0.75\textwidth]{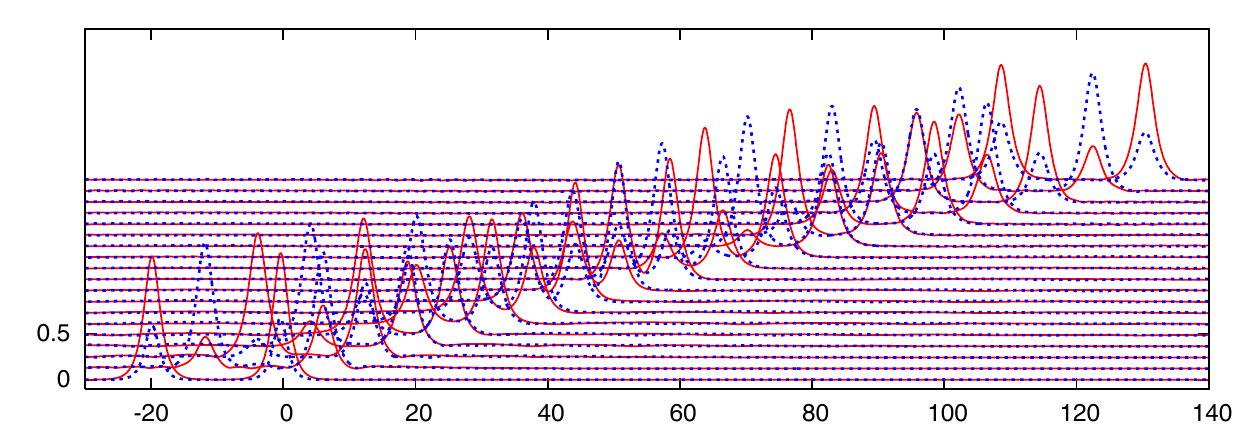}}
\centerline{\small (c) $\deljump=180^\circ$}
\caption{Takeover collision of QPs with initial elliptic polarization: polarization angles $\theta_l=23^\circ 45'$, $\theta_r=25^\circ 52'$, $c_l=1$, $c_r=0.8$, $n_{l\psi}=n_{r\psi}=-1.1$, $n_{l\phi}=n_{r\phi}=-1.5$, $\alpha_1=0.75$, $\Gamma=0.175$, $X_l=-20$, $X_r=0$ and different initial phase differences.}
\label{fig:clashd180ot}
\end{figure}
Recall that these solutions apply specifically to the case $\alpha_2=0$. Our initial condition is valid only for this case. In Figure \ref{fig:clashcird180} we present these features for fixed set of  parameters varying only the initial phase difference $\deljump$. As we discussed in \cite{aip1340} there is no velocity shift on the place of interaction. The individual masses start to oscillate since the same beginning with an equal period and with a phase shift and an amplitude depending on the concrete choice of $\deljump$ The total mass, however, is perfectly conserved -- it is equal to 6.0683. Due to the symmetry the pseudomomentum $P=10^{-3}\div 10^{-11}$ and evidently approximates 0. The result is the same for all considered phase differences ($0^\circ$, $90^\circ$, $180^\circ$) (Figure \ref{fig:clashcird180}). What happens with the energy one can see reviewing the next Figure~\ref{fig:energycird0}. The magnitude of the energy depends essentially on the magnitude of initial phase difference value $\deljump$. The dependence is non-monotonous: $E=-3.078$ for $\deljump=0^\circ$; $E=-4.15$ for $\deljump=90^\circ$; $E=-1.139$ for $\deljump=180^\circ$. Figure~\ref{fig:energycird0} demonstrates cogently the conservation properties of the used difference scheme -- one observes excellent conservation of the total energy of the interacting soliton system during the whole time of evolution. The next Figure~\ref{fig:polarcird180} elucidates the polarization evolution. The individual polarizations oscillate depending on both the magnitude of linear coupling $\Gamma$ and initial phase difference $\deljump$. These oscillations do not depend on the interaction but the period depends on $\Gamma$. The magnitude $\deljump$ influences the amplitude of oscillation. In contrast to masses (see Figure~\ref{fig:massmomcird90}) the total polarization is not constant. It breathes also though with very small amplitude. It is well seeing that on the place of interaction one observes a shock of both the individual and total polarization. The last result can be interpret as a generalization of the polarization conservation law for the case of rotational polarization -- it breathes in the period but it is conserved in the whole period.

\begin{figure}[t!]
\centerline{\includegraphics[width=0.75\textwidth]{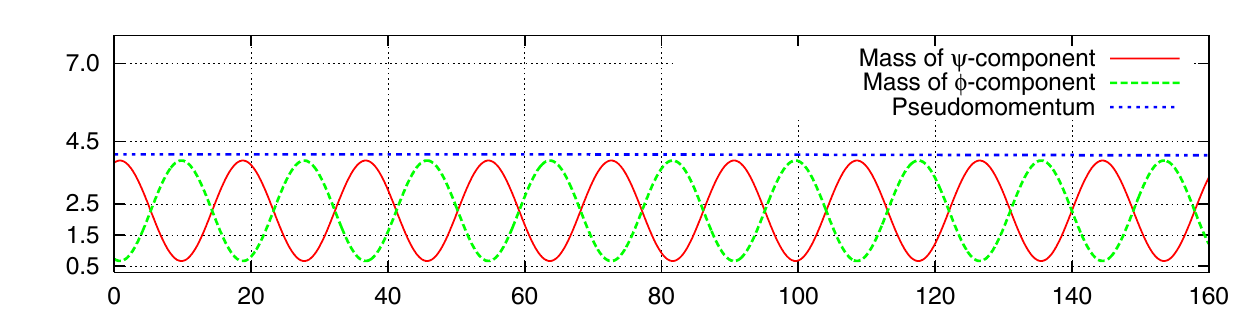}}
\centerline{\small (a) $\deljump=0^\circ$}
\centerline{\includegraphics[width=0.75\textwidth]{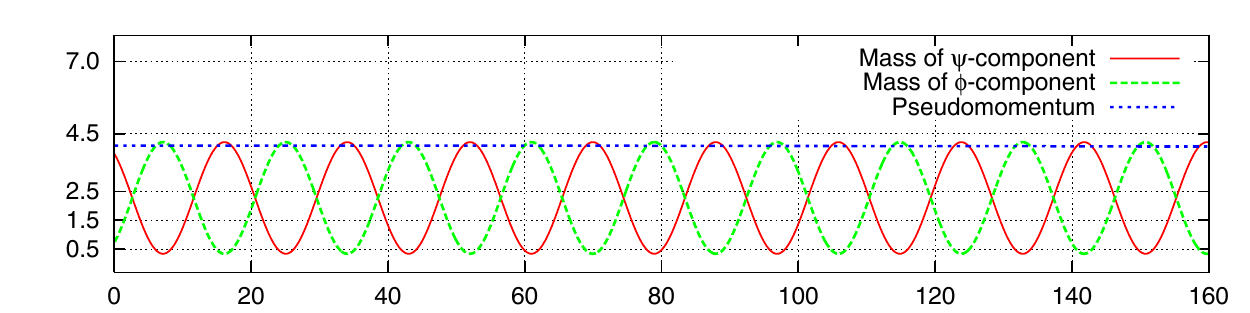}}
\centerline{\small (b) $\deljump=90^\circ$}
\centerline{\includegraphics[width=0.75\textwidth]{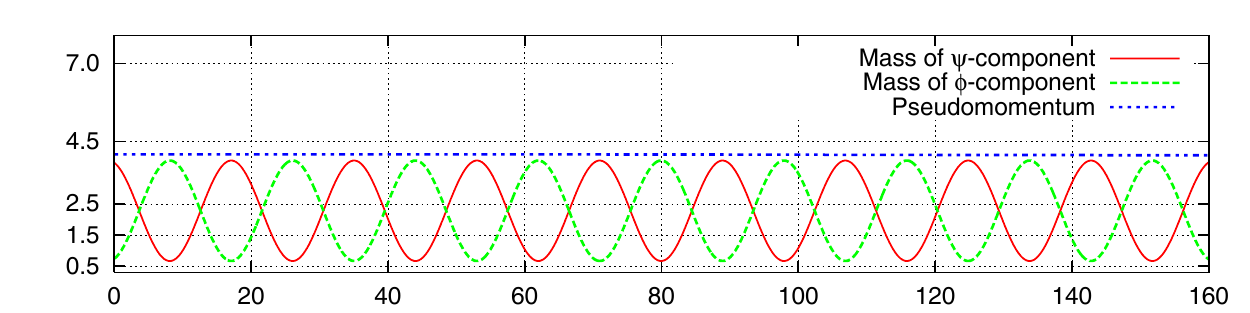}}
\centerline{\small (c) $\deljump=180^\circ$}
\caption{Mass and pseudomomentum dynamics for $\deljump=0^\circ$ corresponding to Figure~\ref{fig:clashd180ot}.}
\label{fig:massmomd180ot}
\end{figure}
\begin{figure}[h!]
\centerline{\includegraphics[width=0.75\textwidth]{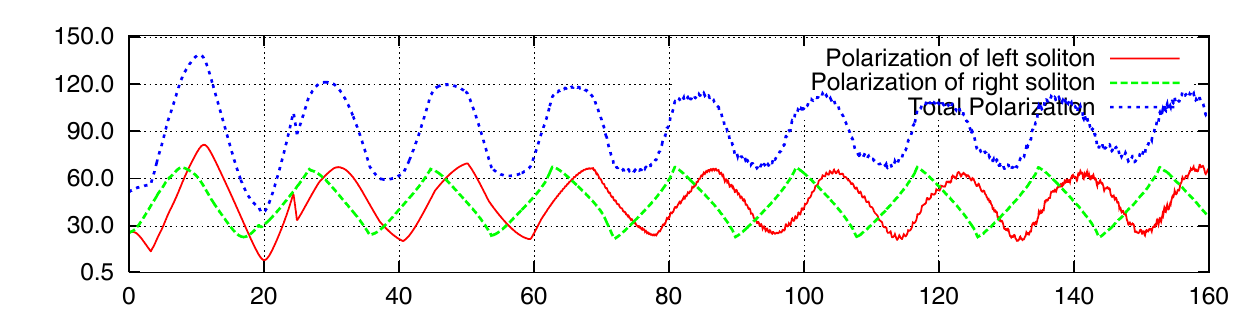}}
\centerline{\small (a) $\deljump=0^\circ$}
\centerline{\includegraphics[width=0.75\textwidth]{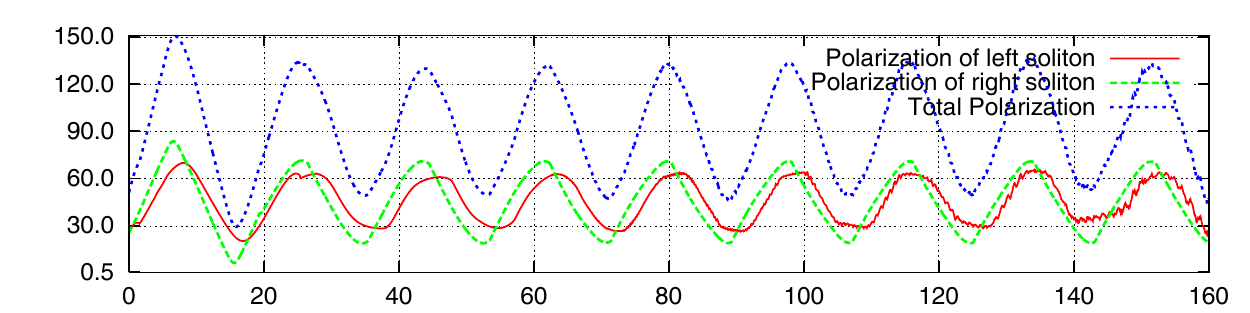}}
\centerline{\small (a) $\deljump=90^\circ$}
\centerline{\includegraphics[width=0.75\textwidth]{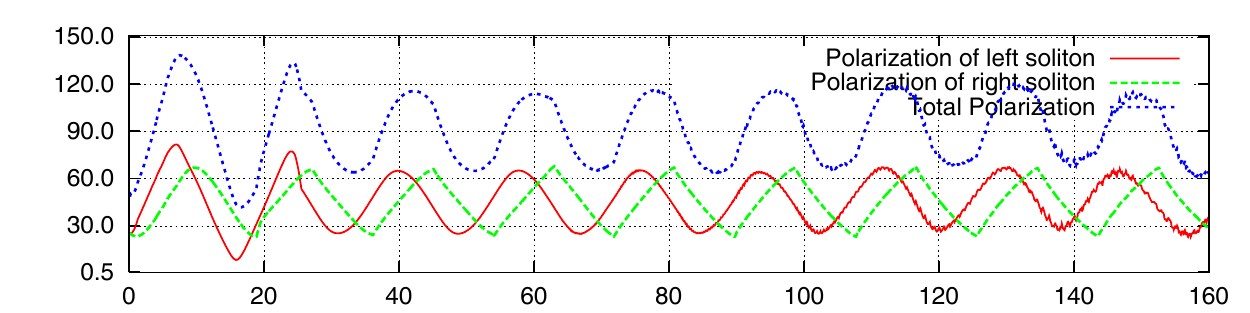}}
\centerline{\small (a) $\deljump=180^\circ$}
\caption{Individual and total polarization dynamics in dependence on initial phase difference corresponding to Figure~\ref{fig:clashd180ot}.}
\label{fig:polard180ot}
\end{figure}
\begin{figure}[h!]
\centerline{\includegraphics[width=0.75\textwidth]{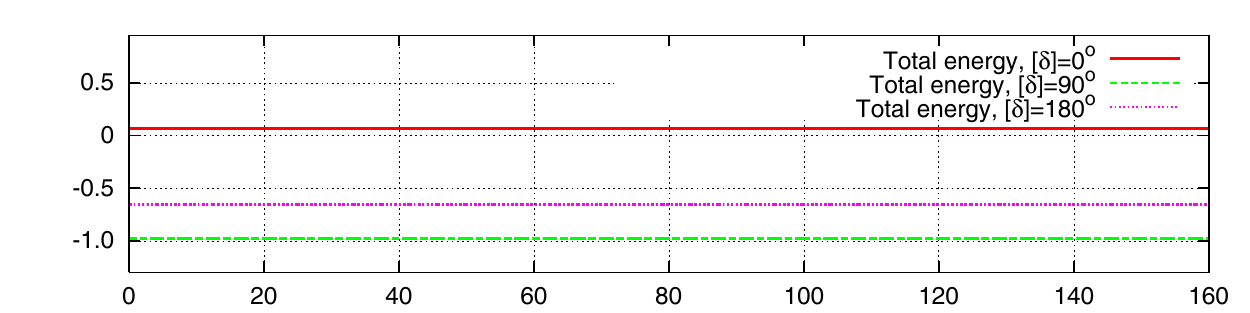}}
\vspace{-0.1in}
\caption{Influence of the initial phase difference on the total energy corresponding to Figure~\ref{fig:clashd180ot}.}
\label{fig:energyd180ot}
\end{figure}

Our observations show that these interactions are independent of initial phase shift and no-phase velocity shift occurs on the place of interaction. The trajectories of QP-centers suffers slight shift. They were investigated and discussed in \cite{aip1340} for \emph{sech}-profiles and initial linear polarization and we do not concern them again.

\subsection{Elliptically polarized solitons. Head-on collision}
        In order to extend our considerations we conduct series of experiments with general initial polarization and nontrivial linear coupling. For convenience we use the same magnitude $\Gamma=0.175$. On Figure \ref{fig:clasheld180} the case of initial elliptic polarization is present for initial phase shifts $\deljump=0^\circ, 90^\circ, 135^\circ, 180^\circ$. Our conclusion is that the linear coupling combined with the initial phase difference $\deljump$ can shift only the soliton trajectories and keeps the phase velocities unchanged. A nontrivial phase shift is possible after involving nonlinear coupling $\alpha_2 \ne 0$ \cite{TodoAIMS10}. 
     We consider two solitons with equal initial elliptic polarization angles of $23^\circ 44'$. The initial configuration is generated from the auxiliary bifurcation system (\ref{eq:conjsys}). The above angle corresponds to the parametric set $n_{l \psi}=n_{r \psi}=-1.1$, $n_{l \phi}=n_{r \phi}=-1.5$, $c_l=-c_r=1$, $\alpha_1=0.75$. Obviously the results are qualitatively the same as in the previous case of circular polarization and they do not require a detailed discussion. On Figure \ref{fig:ell_allenergiies} is plotted the energy depending on the initial phase difference. The correlation between them is essential. The conservation is again excellent: $\deljump=0^\circ$ -- $E=-0.262$; $\deljump=90^\circ$ -- $E=-0.821$; $\deljump=135^\circ$ -- $E=-0.206$; $\deljump=180^\circ$ -- $E=0.640$. The breathing behavior of the total polarization is well visible as well as the affect of the initial phase difference upon the amplitude of the individual polarization angles which result to the magnitude of the total polarization angle. Figure \ref{fig:clasheld180} answers convincingly the question about the dynamics of the rotational polarization when the type of initial polarization is general. The net conservation law is generalized as net breathing (oscillating) law which conserved the total polarization within one time period. This result has not an analytic analogue and enriches the phenomenology of CNLSE. Among the three conservation laws about the mass, pseudomomentum and energy we found after an extensive set of numerical investigations one more net conservation law.

\subsection{Elliptically polarized solitons. Takeover interaction}
The pattern in takeover interaction is the same qualitatively but the required difference in the velocity magnitudes bears bigger trajectory shift of the slower soliton and slower trajectory shift of the bigger soliton (see \cite{aip1340}). Here we consider this kind of interaction in order to implement an interaction of elliptically polarized solitons with different polarization angles. Also such kind of evolution and interaction do not possess a symmetry. Our main goal is to check and confirm the conclusions about the net polarization realized in the previous subsection. We choose $n_{l \psi}=n_{r \psi}=-1.1$, $n_{l \phi}=n_{r \phi}=-1.5$, $c_l=1$, $c_r=0.8$, $\alpha_1=0.75$, $\Gamma=0.175$ and focus on the effects of $\vec \delta$. The above parametric set results in polarization angles $\theta_l=23^\circ 45'$, $\theta_r=25^\circ 52'$ of the initial envelope configuration. The approximate initial condition is generated from the auxiliary bifurcation system. The solitons start with different elliptic polarizations and are an object of rotational polarization as well due to nontrivial linear coupling $\Gamma$ (Figure \ref{fig:clashd180ot}). We find out that the initial phase differences of the components play an essential role on the full energy of QPs (Figure~\ref{fig:energyd180ot}). In particular, when $\deljump=0^\circ$ -- $E=0.0673$; $\deljump=90^\circ$ -- $E=-0.976$; $\deljump=180^\circ$ -- $E=-0.657$. Just in opposite, the pseudomomentum is also conserved and does not depend on the initial phase difference. It is not trivial because of the nonsymmetry of the QP configuration and its net values  $P=4.08\div4.05$ (Figure~\ref{fig:massmomd180ot}). The individual masses, however, breathe together with the individual (rotational) polarizations. Their amplitude and period do not influenced from the initial phase difference (Figure~\ref{fig:massmomd180ot}) and are conserved within one full period of the breathing. The total mass is constant. Both the individual and total polarizations breathe and suffer a `shock in polarization' when QPs enter the collision. The polarization amplitude evidently depends on the initial phase difference (Figure~\ref{fig:polard180ot}). The above quantities are conserved within one full period of the breathing. In all considered cases of evolution and interaction the polarization angle of QPs can change independently of the collision due to the real linear coupling $\Gamma$.

\section{Conclusion}
 In this work we are aiming to complete the information about the dynamical properties of the total polarization influenced by the initial parameters and, in particular, by the initial polarization of general type as well as by the nontrivial linear coupling $\Gamma$ and the initial phase difference $\deljump$. For the case of general elliptic polarization there is not available an exact initial condition concerning the shape of solitary envelopes. Solving an auxiliary bifurcation system of ordinary differential equations we construct numerically an approximate initial conditions for the full range of polarization angles between $0^\circ$ and $90^\circ$ and with non-$sech$ shape. Let us denote that the particular cases $\theta=0^\circ;90^\circ$ (left and right linear polarization) when one has an exact $sech$-like initial condition were implemented in detail in \cite{aip1340}. We consider the two main interactions -- head-on and takeover collisions of envelopes  when only a real linear coupling $\Gamma$ is present (the cross-modulation parameter $\alpha_2=0$). Adding an imaginary part in $\Gamma$ leads to gain/blow-up of the solution accompanied by a violation of the conservation laws (see \cite{aip1340}) and we will not concern again this property in the present work. The main conclusions can be summarized as follows
 \begin{itemize}
 \item The interaction conserves perfectly the pseudomomenta and energy. The pseudomomenta are not influenced by the initial phase difference  $\deljump$ while the magnitude of the energy is strongly dependent on it. The variation is nonmonotonous.

 \item The total mass is also perfectly conserved but this is realized on the inner share and exchange of local masses of both components of QP periodically in the time. This periodicity is governed by the real part of the linear coupling $\Gamma$ and generates breathing of the masses among with the breathing of the solitary envelopes. The mass breathing does not depend on the interaction and is present during the whole time of evolution. The amplitude of mass oscillation depends on the initial phase difference.

 \item  The individual polarizations breathe similar to individual masses but  in opposite to them with a phase shift. This is the reason  the total rotational polarization to be conserved within one period. On the place of interaction the individual and total polarization angles suffer a discontinuity keeping the else functions smooth. In other words the total polarization breathes also and a bifurcation occurs during the cross-section of the interaction. And this is the main result of this investigation. The amplitudes of the individual polarizations are affected by the initial phase difference.

 \item We found and generalized for the case of rotational polarization one more conservation law -- for the net total polarization. So, the linearly CNLSE possess four conservation laws: for the total mass, for the pseudomomentum, for the energy, and for the magnitude of the total polarization. Let us emphasize that the last law has not an analytic analogue and it unearths more intrinsic information for the propagating and interacting solitons as quasi-particles.
 \end{itemize}

 \section{Acknowledgements}
 This investigation is partially supported by the Science Fund of Ministry of Education, Science and  Youth of Republic Bulgaria under grant DDVU02/71.

 These results were reported on the symposium ``Recent Trends in Nonlinear Wave Phenomena: Achievements and Challenges'' within the 6th IMACS International Conference on  Nonlinear Evolution Equations and Wave Phenomena: Computation and Theory, Athens, GA, April 04-07, 2011. The paper is devoted to my mentor and supervisor Professor Christo I. Christov on occasion of his 60th birthday.

\bibliographystyle{model1b-num-names}


\begin{thebibliography}{10}
\normalsize
\bibitem{2}
K.~W. Chow, K.~Nakkeeran, and B.~Malomed.
\newblock Periodic waves in bimodal optical fibers.
\newblock {\em Opt.Comm.~}, 219:251--9, 2003.

\bibitem{Chri_Reading}
C.~I. Christov.
\newblock {\em Gaussian elimination with pivoting for multidiagonal systems,
  Internal Report {\bf 4}}.
\newblock University of Reading, 1994.

\bibitem{ChriDostMau}
C.~I. Christov, S.~Dost, and G.~A. Maugin.
\newblock Inelasticity of soliton collisions in system of coupled {NLS}
  equations.
\newblock {\em Physica Scripta}, 50:449--454, 1994.

\bibitem{collings}
B.~C. Collings, S.~T. Cundiff, N.~N. Akhmediev, J.~M. Soto-Crespo, K.~Bergman,
  and W.~H. Knox.
\newblock Polarization-locked temporal vector solitons in a fiber laser:
  experiment.
\newblock {\em J. Opt. Soc. Am. B}, 17:354--365, 2000.

\bibitem{hempelmann}
U.~Hempelmann.
\newblock Polarization coupling and transverse interaction of spatial optical
  solitons in a slab waveguide.
\newblock {\em J. Opt. Soc. Am. B}, 12:77--86, 1995.

\bibitem{lakoba}
T.~I. Lakoba, D.~J. Kaup, and B.~A. Malomed.
\newblock Solitons in nonlinear fiber couplers with two orthogonal
  polarizations.
\newblock {\em Phys. Rev. E}, 55:6107--6120, 1997.

\bibitem{Optical_IMACS}
W.~J. Sonnier and C.~I. Christov.
\newblock Strong coupling of {Schr\"odinger} equations: Conservative scheme
  approach.
\newblock {\em Math. Comp. Simul.}, 69:514--525, 2005.

\bibitem{TahaAblovitz}
T.~R. Taha and M.~J. Ablovitz.
\newblock Analytical and numerical aspects of certain nonlinear evolution
  equations. {II}. {N}umerical, {S}chr\"odinger equation.
\newblock {\em J. Comp. Phys.}, 55:203--230, 1984.

\bibitem{TodoAIMS10}
M.~D. Todorov.
\newblock Polarization dynamics during takeover collisions of solitons in
  systems of {C}oupled {N}onlinear schr{\"o}dinger {E}quations.
\newblock {\em Discrete and Continuous Dynamical Systems}, Supplement, 2011.
\newblock in press.

\bibitem{TodoChri06}
M.~D. Todorov and C.~I. Christov.
\newblock Conservative scheme in complex arithmetic for {C}oupled {N}onlinear
  {S}chr\"odinger {E}quations.
\newblock {\em Discrete and Continuous Dynamical Systems}, Supplement:982--992,
  2007.

\bibitem{TodoChri08}
M.~D. Todorov and C.~I. Christov.
\newblock Collision dynamics of circularly polarized solitons in nonintegrable
  {C}oupled {N}onlinear {S}chr{\"o}dinger {S}ystem.
\newblock {\em Discrete and Continuous Dynamical Systems}, Supplement:780--789,
  2009.

\bibitem{bgsiam09}
M.~D. Todorov and C.~I. Christov.
\newblock On the solution of the system of {ODE}s governing the polarized
  stationary solutions of {CNLSE}.
\newblock In {\em 3rd Annual Session of BGSIAM, December 22-23, 2008}, pages
  83--86. Demetra Ltd., 2009.

\bibitem{matcom10}
M.~D. Todorov and C.~I. Christov.
\newblock Collision dynamics of elliptically polarized solitons in {C}oupled
  {N}onlinear {S}chr{\"o}dinger {E}quations.
\newblock {\em Math. Comput. Simul.}, doi:10.1016/j.matcom.2010.04.022, 2010.

\bibitem{aip1340}
M.~D. Todorov and C.~I. Christov.
\newblock Collision dynamics of polarized solitons in linearly {CNLSE}.
\newblock In {\em Intl. Workshop on Complex Structures, Integrability and
  Vector Fields, AIP CP1340}, pages 144--153. American Institute of Physics,
  Melville, NY, doi: 10.1063/1.3567133, 2011.

\bibitem{1}
H.~Winful.
\newblock Self-induced polarization changes in biregringent optical fibers.
\newblock {\em Appl.Phys.Lett.~}, 47(3):213--5, 1985.

\bibitem{yang}
J.~Yang.
\newblock Classification of the solitary waves in {C}oupled {N}onlinear
  {S}chr{\"o}dinger {E}quations.
\newblock {\em Physica D}, 108:92--112, 1997.

\end{thebibliography}

\end{document}